\documentclass[10pt]{article}
\usepackage{amsmath,amssymb,amsfonts,amsthm}
\usepackage[matrix,arrow,curve]{xy}
\author{Victor Tourtchine\thanks{Partially supported by the
grants NSH-1972.2003.01}\\
{\footnotesize\it Independent University of Moscow}\\
{\footnotesize\it Universit\'e Catholique de Louvain}\\
{\footnotesize\tt turchin@math.ucl.ac.be, vitia-t@yandex.ru}}
\title{What is one-term relation for higher homology of long knots}
\date{}

\newtheorem*{theorem*51}{Theorem \ref{t51}}
\newtheorem{theorem}{Theorem}[section]
\newtheorem{lemma}[theorem]{Lemma}

\newtheorem{remark}[theorem]{Remark}
\newtheorem{remark-corollary}[theorem]{Remark-Corollary}
\newtheorem{definition}[theorem]{Definition}
\newtheorem{corollary}[theorem]{Corollary}
\newtheorem{proposition}[theorem]{Proposition}

\newtheorem{notation}[theorem]{Notation}

\def\mod{\mathop{\rm mod}\nolimits}
\newcommand\R{{\mathbb R}}

\newcommand\Z{{\mathbb Z}}
\newcommand\Q{{\mathbb Q}}
\newcommand\N{{\mathbb N}}

\newcommand\ZZ{{\mathcal Z}}

\newcommand\BB{{\mathcal B}}
\newcommand\WW{{\mathcal W}}

\newcommand\TD{{\mathcal {TD}}}
\newcommand\TstD{{\mathcal {T}_*\mathcal{D}}}
\newcommand\TzD{{\mathcal {T}_0\mathcal{D}}}
\newcommand\TdstD{{\mathcal T}^*_{\,*}{\mathcal D}}

\newcommand\ie{{\it i.e. }}
\newcommand\cf{{\it cf.}}

\newcommand\DTheta{\unitlength=0.13em
\begin{picture}(10,10)
 \put(5,5){\circle{10}}
 \put(0,5){\line(1,0){10}}
\end{picture}}

\def\id{\mathop{\rm id}\nolimits}

\def\nbox{\quad$\Box$}
\def\nboxm{\quad\Box}

\newcommand\rth{\refstepcounter{equation}}
\newcommand\numb{\rth{\rm \theequation}}
\numberwithin{equation}{section}

\begin{document}
\maketitle \sloppy

{\footnotesize
\begin{abstract}
Vassiliev's spectral sequence for long knots is discussed. Briefly speaking we study what happens if the strata of non-immersions are ignored. 

Various algebraic structures on the spectral sequence are introduced. General theorems about  these structures imply, for example, that the bialgebra of chord diagrams is 
polynomial for any field of coefficients.
\end{abstract}


\medskip

}

\subsection*{Acknowledgement}
The author thanks V.~Vassiliev, P.~Lambrechts, M.~Liverenet, F.~Patras, B.~Fresse, S.~Lando for interesting discussions. The author is greatful to the Universit\'e Catholique de Louvain
where this paper was written for hospitality.

\section{Introduction}\label{introduction}
\subsection{Bialgebra of chord diagrams. Knots and framed knots}\label{ss11}
Bialgebra of chord diagrams $\BB_{0}$ is the dual space to the space $\WW_{0}$ of weight systemes associated to finite type knot invariants, \cf~\cite{BarNatan}.
$\BB_{0}$ is described as a space spanned  by chord diagrams on a circle and quotiented out by 4-term relations ($4T$) and by one-term relations ($1T$). It is also natural to 
consider a similar bialgebra $\BB$ which is obtained by taking a quotient only by $4T$-relations. Bialgebra $\BB$ is  the dual space to the space $\WW$ of weight systems associated to
finite type invariants of framed knots\footnote{For example, weight systems associated to representations of compact Lie algebras respect only $4T$-relations. The corresponding quantum
invariants are defined on framed knots, \cf~\cite{Kassel}.}. These bialgebras are almost the same. One has:
$$
\BB\simeq\BB_{0}\otimes\Z[\Theta]  
\eqno(\numb)\label{eq11}
$$
where $\Z[\Theta]$ is a polynomial bialgebra with the only primitive generator $\Theta$ which is the only 1-chord diagram \DTheta.\footnote{Over $\Q$ this is proven 
in~\cite{BarNatan}. Over $\Z$ this is a consequence 
of Proposition~6.3.8 in~\cite{LZ}. Splitting~\eqref{eq11} is also a consequence of Theorem~\ref{t51}, see below.}

In this paper we will consider spaces $Emb_{d}$, $d\geq 3$, (sometimes we will write  simply $Emb$) of {\it long knots}, \ie
spaces of smooth embeddings $\R^1\hookrightarrow\R^d$ coinciding with a fixed linear map outside some compact set of $\R^1$. For $d\geq 4$, this space is connected 
(and even simply connected). The connected components of $Emb_{3}$ define isotopy classes of long knots which are in one-to-one correspondence with the isotopy classes of usual
knots. Analogously the bialgebra $\BB_{0}$ (or $\BB$) of chord diagrams on a circle and the corresponding bialgebra on a line coincide. This fact is a consequence of $4T$-relations.

\begin{figure}[!h]
$$
  \begin{picture}(50,15)
    \put(0,0){\line(1,0){50}}
    \qbezier(5,0)(20,15)(35,0)
    \qbezier(15,0)(30,15)(45,0)
    \qbezier(25,0)(35,11)(45,0)
  \end{picture}
 = 
  \begin{picture}(50,15)
    \put(0,0){\line(1,0){50}}
    \qbezier(5,0)(20,15)(35,0)
    \qbezier(15,0)(25,11)(35,0)
    \qbezier(25,0)(35,11)(45,0)
  \end{picture}
\, \mod 4T.
 $$
\caption{Example of a chord diagram on the line}\label{fig1}
\end{figure}

For the spaces $Emb_{d}$, $d\geq 4$, the bialgebra $\BB_{0}$ is a subbialgebra of the homology bialgebra $H_{*}(Emb_{d})$, \cf~\cite{CCL,CCL2}.
\footnote{To be precise in the case of even $d$, bialgebra $\BB_{0}$ must be replaced by some its non-trivial super-analogue. 
Bialgebra structure on $H_{*}(Emb_{d})$ is induced by $H$-space structure of $Emb_{d}$, multiplication of knots being concatenation.}

\subsection{Approach of V.~Vassiliev. Ignoring strata of non-immersions}\label{ss12}
Historically, the finite type knot invariants, weight systems and the bialgebra of chord diagrams appeared in the approach of V.~Vassiliev of studying the spaces of knots, \cf~\cite{V1,V4}.
Following V.~Vassiliev, space $Emb_{d}$ is an open everywhere dense subset in the infinite-dimensional affine space $\R^{\omega d}$ of all smooth maps $\R^1\to\R^d$ coinciding
with a fixed linear map outside some compact set of $\R^1$. The complement $\Sigma_{d}=\R^{\omega d}\setminus Emb_{d}$ of this open subset is called the {\it discriminant space}
or simply {\it discriminant}. It consists of maps having self-intersections and/or singularities. The idea of Vassiliev is to use the Alexander duality:
$$
\tilde H^*(Emb_{d})\equiv \tilde H^*(\R^{\omega d}\setminus\Sigma_{d})\simeq\tilde H_{\omega d-*-1}(\bar\Sigma_{d}),
$$
and express cohomogy classes of $Emb_{d}$ as linking number with appropriate chains (of finite codimension) lying in the one-point compactification
$\bar\Sigma_{d}$ of the disriminant~\footnote{We refer the reader to~\cite{V1,V4} for rigorous statements.}. 

The main tool of V.~Vassiliev is a spectral sequence associated to a natural filtration of resolution of $\Sigma_{d}$. For $d\geq 4$, this spectral sequence
does compute the cohomology groups of $Emb_{d}$.

\begin{notation}\label{n12}
{\rm In this paper we will consider several spaces which are also complements to discriminants in some affine functional spaces. The spectral sequence obtained
by the above method will be called {\it Vassiliev spectral sequence} associated to the corresponding space. \nbox}
\end{notation}

In Vassiliev's approach, chord diagrams correspond to strata of maps $f\colon\R^1\to\R^d$ in $\Sigma_{d}$ having a finite number of double self-intersections.
$4T$-relations and $1T$-relations correspond respectively to triple self-intersections and degeneration of the first differential. So, the rejection of the $1T$-relations 
can be viewed as ignoring of strata of non-immersions in the discriminant. But note that the strata of non-immersions lie in the closure of the strata of immersions. 
So, we speak only about some formal ignoring of the above strata in our algebraic spectral sequence calculations.

In~\cite{T5} we described complex $\TstD$ (sometimes we will denote it by $\TstD^d$ when we will need to emphasize the ambient dimension $d$) which
is quasi-isomorphic to the $E_{0}$ term of the Vassiliev spectral sequence and serves to simplify the calculation of the first term. 
One can recover this complex even in~\cite{V1}, but it is not defined explicitly there. One is referred to use the results of~\cite{V3} to construct this complex.
In~\cite{T5} we defined also a quotient-complex $\TD$ of $\TstD$ spanned by the diagrams corresponding to strata of immersions.

The main aim of this paper is to describe how the homology of $\TstD$ is related to that of $\TD$.

\begin{theorem*51}
Complex $\TD$ is quasi-isomorphic to the tensor product $\TstD\otimes\ZZ$, where $\ZZ$ is a  complex computing the cohomology of 
$\Omega^2S^{d-1}$.\nbox
\end{theorem*51}

Bialgebra  $\BB$ (resp. $\BB_{0}$), is a subspace of the homology of the complex dual to $\TD^{odd}$ (resp. $\TstD^{odd}$).
So, Theorem~\ref{t51} is a non-trivial higher (co)homology generalization of~\eqref{eq11}. Indeed, splitting~\eqref{eq11} is a consequence of the Kunneth formula.
Note, that the $Tor$-part of the Kunneth formula is absent in the bigradings of $\BB$ (resp. $\BB_{0}$).
In the above theorem $\Theta$ is interpreted as a generator of 
$H_{d-3}(\Omega^2 S^{d-1})\simeq\Z$.

In the case when $d$ is even, Theorem~\ref{t51} implies an analogous splitting:
$$
\tilde\BB\simeq\tilde\BB_{0}\otimes \left(\Z[\Theta]/[2\Theta^2{=}0]\right).
$$
$\tilde\BB$ (resp. $\tilde\BB_{0}$) designates the super-analogue of $\BB$ (resp. $\BB_{0}$), \cf~\cite[Section~35.2]{T2}, \cite[Section~6]{T4}. 
Note that $\Theta$ is now of odd degree, and 
$\Z[\Theta]/[2\Theta^2{=}0]$ is  the free super-commutative algebra over $\Z$ of one odd generator.

\subsection{Main results}\label{ss13}
As we mentioned in the previous subsection, Theorem~\ref{t51} is the main result of the paper. Other results are given by Theorems~\ref{t93}, \ref{t95}, \ref{t113} and 
Corollaries~\ref{c114}, \ref{c115}. Theorems~\ref{t93},~\ref{t95} are given without proof and present results of calculations. We give there explicit combinatorial
formula in the spirit of~\cite{V5,V6} for some non-trivial cocycles of the space $Emb^+$. Space $Emb^+$ is defined in the next section and is 
a \lq\lq good replacement" for the space of framed knots.
The results of Section~\ref{s11} show some freeness proprieties of the first term of the Vassiliev spectral sequence for long knots. As a consequencethe the  bialgebra of
chord diagrams is a free polynomial bialgebra for any field of coefficients (Corollary~\ref{c115}).

\section{Homotopy fiber $Emb_{d}^+$}\label{s2}
As we have mentioned in Section~\ref{introduction}, the rejection of the $1T$-relations (ignoring strata of non-immersions) was always
associated with the space $Emb_{3}^{framed}$ of framed knots in $R^3$. It turns out that this association is wrong when we want to consider 
higher (co)homology of the space $Emb_{3}$ (or more generally of $Emb_{d}$). The space that should replace $Emb_{d}^{framed}$ is the homotopy fiber
$Emb_{d}^+$ of the inclusion $Emb_{d}\hookrightarrow Imm_{d}$, where $Imm_{d}$ designates the space of long immersions, \ie immersions 
$\R^1\to\R^d$ with the fixed linear behavior at infinity.

The idea to consider the space $Emb^+$ appeared first in~\cite{Vol,Sinha2}, where the authors studied Goodwillie's approach of calculus of embeddings, see also Section~\ref{s10}. 

\begin{theorem}\label{t21}
{\rm D.~Sinha~\cite{Sinha2}.} The inclusion $Emb_{d}\hookrightarrow Imm_{d}$, $d\geq 3$, is a homotopy trivial map, and therefore
its homotopy fiber $Emb_{d}^+$ is weakly homotopy equivalent to the direct product
$$
Emb_{d}^+\simeq Emb_{d}\times\Omega Imm_{d}\simeq Emb_{d}\times\Omega^2 S^{d-1}. \nboxm 
\eqno(\numb)\label{eq22}
$$ 
\end{theorem}

\noindent {\bf Proof:} The proof of D.~Sinha is so simple that we repeat it. $Imm_{d}$ is naturally mapped to $\Omega S^{d-1}$: it is defined by taking the direction
of the derivative. It is an easy exercise to construct a homotopy inverse to this map, and hence $Imm_{d}$ is homotopy equivalent to $\Omega S^{d-1}$.
Consider the composite map $Emb_{d}\to\Omega S^{d-1}$. Let us show that is contractible. 
The contraction is defined as follows (the parameter $a\in[0,+\infty]$ defines the homotopy):
$$
\frac{f(t+a)-f(t-a)}{||f(t+a)-f(t-a)||},
$$
where $f\in Emb_{d}$. \nbox

\bigskip

We will see below that Theorem~\ref{t51} is a manifistation of the homotopy splitting~\eqref{eq22}.

It is well known that $Emb_{3}^{framed}$ is homotopy equivalent to $Emb_{3}\times\Z$. So, from Theorem~\ref{t21} we see that 
the space $Emb_{3}^+$ is a good replacement for $Emb_{3}^{framed}$ on the level of isotopy classes:
$\pi_{0}(Emb_{3}^+)=\pi_{0}(Emb_{3})\times\Z=\pi_{0}(Emb_{3}^{framed})$.

\section{$Emb_{d}^+$ from the point of view of Vassiliev theory}\label{s3}
Space $Emb_{d}^+$ can be regarded as a space of maps
$$
g\colon [0,1]\times\R\to\R^d,
$$
such that $g(0,-)\colon\R^1\to\R^d$ is an embedding (in $Emb_{d}$). $g(1,-)$ is the fixed linear map, and $g(a,-)$ is an immersion (in $Imm_{d}$)
for any $a\in (0,1)$. Space $Emb_{d}^+$ is a Serre fibration over $Emb_{d}$. Theorem~\ref{t21} asserts that this fibration is homotopy trivial.

Space $Emb_{d}^+$ is an open everywhere dense subset in the affine space of all maps $g\colon [0,1]\times\R\to\R^d$ verifying

1) $g(1,-)$ is the fixed linear map;

2) $g(a,t)=g(1,t)$ outside some compact subset of $[0,1]\times\R$.

\smallskip

The complement $\Sigma_{d}^+$ consists of maps such that $g(0,-)$ is not an embedding or $g(a,-)$ is not an immersion for some $a\in (0,1)$.

\medskip

The aim of this section is to define a new complex $\TdstD$ which is quasi-isomorphic to $\TD$.

\begin{proposition}\label{p31}
The defined below complex $\TdstD$ is quasi-isomorphic to the $E_{0}$ term of the Vassiliev cohomological spectral sequence associated 
to the space $Emb_{d}^+$ (see Notation~\ref{n12}).\nbox
\end{proposition}

\noindent {\bf Sketch of the proof:} This result is a standard application of Vassiliev's theory of discriminants, \cf~\cite{V4}. To be precise complexe $\TdstD$ is obtained as follows:

First, one considers simplicial resolution $\sigma_{d}^+$ of $\Sigma_{d}^+$. Second, $\sigma_{d}^+$ has a natural filtration,
that defines a spectral sequence (Vassiliev spectral sequence). Third, we introduce an auxiliary filtration in the terms of the main filtration 
in order to compute the first term of the main spectral sequence. The auxiliary spectral sequence asociated to the auxiliary filtration collapses in the second term.
Its first term is exactly our complex $\TdstD$. \nbox

\bigskip

The space of the complex $\TdstD$ is spanned by {\it $T_{*}^*$-diagrams} and quotiented out by 3-term relations (these relations are the Arnol'd's relations appearing in the
cohomology of configuration  spaces, \cf~\cite{A1}).

\begin{definition}\label{tdstd}{\rm
Any {\it $T_*^*$-diagram} is the following set of data:

1) $n\geq 0$  points $t_{1}<t_2<\ldots <t_n$ on the line $\R^1$. We will call them {\it active points} of the diagram.

2) A set of oriented chords joining these points. Each chord joins two distinct active point of the diagram. Two distinct points are joined  by no more then one chord
(no multiple edges).

3) From each active point we draw up a vertical half-line. On each half-line we fix some number of points and put an asterisk in each of them. The
asterisks can coincide with the active points of the diagram or can be over them. The asterisks of the first type will be called {\it bottom asterisks};
the asterisks of the second type will be called {\it top asterisks}.

4) We demand that each active point is adjacent to a chord, or its half-line contains at least one asterisk.

5) We demand that if we remove the line $\R^1$, all the half-lines and all the asterisks, the remaining graph would be a forest, \ie a disjoint union of trees. \nbox}
\end{definition}

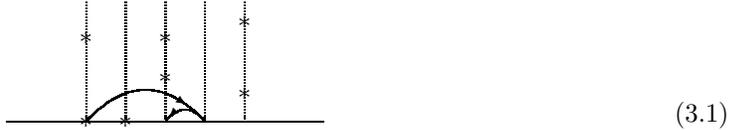
$$
\unitlength=0.3em
\begin{picture}(40,15)
\put(0,0){\line(1,0){40}}
\multiput(10,0)(5,0){5}{\qbezier[40](0,0)(0,7.5)(0,15)}
\put(9.2,-0.75){$*$}
\put(9.2,9.75){$*$}
\put(14.2,-0.75){$*$}
\put(19.2,4.75){$*$}
\put(19.2,9.75){$*$}
\put(29.2,2.75){$*$}
\put(29.2,11.75){$*$}
\qbezier(10,0)(17.5,8)(25,0)
\qbezier(25,0)(22.5,3)(20,0)
\put(22.8,2.1){\vector(3,-2){0}}
\put(20.5,0.5){\vector(-1,-1){0}}
\end{picture}
\eqno(\numb)\label{eq31}
$$ 
\vspace{0.3cm}

Later on we will also need to consider {\it generalized $T_{*}^*$-diagrams}.
\begin{definition}\label{gtdstd}
{\rm
A {\it generalized $T_{*}^*$-diagrams} is the same set of data as a $T_{*}^*$-daigram except that we demand no more Condition 4). In other words we permit active points 
without asterisks and without adjacent chords. \nbox
}
\end{definition}

Each digram is {\it oriented} by ordering of its {\it orienting set}. The orienting set consists of the following objects:

1) Active points $t_{1},\ldots,t_{n}$.

2) Points $a_{i}^j$, $1\leq i\leq n$, on half-lines containing top asterisks (usually we count them from top to bottom).

3) Elements $\alpha_{ij}$ corresponding to oriented chords joining $t_{i}$ with $t_{j}$.

4) Elements $\alpha_{i}^*$ corresponding to bottom asterisks.

5) Elements $\alpha^*{}_{i}^j$ corresponding to top asterisks.

The elements 1)-2) are of degree $-1$, elements 3)-5) are of degree $d-1$. Orientation of a diagram is a monomial including
all the elements 1)-5). For example, orientation of~\eqref{eq31} can be given by the following monomial:
$$
t_1t_2t_3t_4t_5a_{31}a_{32}a_{51}a_{52}\alpha_{14}\alpha_{43}\alpha_{1}^*\alpha_{2}^*\alpha^*{}_{1}^1
\alpha^*{}_{3}^1\alpha^*{}_{3}^2\alpha^*{}_{4}^1\alpha^*{}_{4}^2.
$$
Changing of order of the elements in the orienting monomial is equivalent to multiplication by $\pm 1$ according to the standard
graded commutativity rule.

Space of the complex $\TdstD$ is spanned by the above diagrams and quotiented out by the relations:
$$
\alpha_{ij}=(-1)^{d}\alpha_{ji},
$$
changing of orientation of a chord is equivalent to multiplication by $(-1)^d$; and by the (Arnol'd) relations:
$$
\alpha_{ij}\alpha_{jk}+\alpha_{jk}\alpha_{ki}+\alpha_{ki}\alpha_{ij}=0.
$$
This relation means that the sum of three diagrams which are almost the same except two chords (and these two chords are respectively
$\alpha_{ij}$, $\alpha_{jk}$ for the first diagram; $\alpha_{jk}$, $\alpha_{ki}$ for the second one; and $\alpha_{ki}$, $\alpha_{ij}$ for the third one)
is zero.

The space of $T_{*}^*$-diagrams is bigraded: by {\it complexity} $i$ --- total number of chords and asterisks (elements 1)-3)), and by the {\it number $j$ of geometrically
distinct points} --- total number of active points and top asterisks (elements 1)-2)). Complexity $i$ is the \lq\lq complexity" of the corresponding strata of $\Sigma_{d}^+$.
Number $j$ is the number of the \lq\lq degree of freedom" of the diagram. Bigrading $(p,q)$ of the Vassiliev spectral sequence is related to this bigrading as follows:
$$
p=-i, \quad q=id-j. 
\eqno(\numb)\label{eq32}
$$ 
The corresponding cohomology degree $p+q=i(d-1)-j$, \ie the total degree of the orienting monomial.

We will see that differential $\partial$ of $\TdstD$ conserves the first grading $i$ and diminishes by one the second grading $j$. Differential $\partial$ is a sum
$$
\partial=\partial_{h}+\partial_{v},
$$
where $\partial_{v}$ corresponds to vertical gluings, $\partial_{h}$ --- to horizontal gluings.

The part $\partial_{v}$ is the sum over all top asterisks near to the line $\R^1$, which approach this line and become bottom asterisks:
$$
\unitlength=0.25em
 \partial_{v}\Bigl(\,
 \begin{picture}(20,12)(0,3)
   \put(0,0){\line(1,0){20}}
   \multiput(5,0)(5,0){3}{\qbezier[30](0,0)(0,6)(0,12)}
   \put(4,2.8){$*$}
   \put(4,6.7){$*$} 
   \put(9,4.7){$*$}
   \put(14,4.9){$*$}  
   \put(14,-1){$*$}
   \qbezier(5,0)(7.5,3)(10,0)
   \put(9.5,1){\vector(3,-1){0}}       
 \end{picture}
\, \Bigr)=
\pm\begin{picture}(20,12)(0,3)
   \put(0,0){\line(1,0){20}}
   \multiput(5,0)(5,0){3}{\qbezier[30](0,0)(0,6)(0,12)}
   \put(4,-1){$*$}
   \put(4,6.7){$*$} 
   \put(9,4.7){$*$}
   \put(14,4.9){$*$}  
   \put(14,-1){$*$}
   \qbezier(5,0)(7.5,3)(10,0)
   \put(9.5,1){\vector(3,-1){0}}       
 \end{picture}
 \pm
 \begin{picture}(20,12)(0,3)
   \put(0,0){\line(1,0){20}}
   \multiput(5,0)(5,0){3}{\qbezier[30](0,0)(0,6)(0,12)}
   \put(4,2.8){$*$}
   \put(4,6.7){$*$} 
   \put(9,-1){$*$}
   \put(14,5){$*$}  
   \put(14,-1){$*$}
   \qbezier(5,0)(7.5,3)(10,0)
   \put(9.5,1){\vector(3,-1){0}}       
 \end{picture}
$$
If an active point contained already a bottom asterisk, such a gluing gives zero. Orientation of a diagram of the border is obtained by putting the element
$a_{i}^j$ (corresponding to the asterisk approaching the line) to the first place in the orientation monomial (this gives a sign) and then removing it.

The horizontal part $\partial_{h}$ of the differential is the sum over all possible gluings of couples of neighbor active points $(t_{i},t_{i+1})$ in the diagram. 
If $t_{i}$ and $t_{i+1}$ where joined by a chord $\alpha_{i,i+1}$ (resp. $\alpha_{i+1,i}$),  then this chord becomes bottom asterisk $\alpha_{i}^*$ (resp. $(-1)^d\alpha_{i}$)
unless $t_{i}$ or $t_{i+1}$ contained already an asterisk. In the latter case the gluing gives zero. A gluing gives zero also in two cases: 1)~both points $t_{i}$ and $t_{i+1}$ contain 
asterisks; 2)~the case when $t_{i}$ and $t_{i+1}$ are not joined by a chord
but by a sequence of chords, since such a gluing provides a graph whose one connected component is no more a tree. 

If $t_{i}$ and $t_{i+1}$ contained respectively $k$ and $n$ top asterisks, then after the gluing the half-lines collide, and the result is a sum over ${k+n}\choose{k}$ possible
shuffles of these top asterisks:
$$
\unitlength=0.22em
\partial_h \Bigl(\,
\begin{picture}(15,12)(0,3)
\put(0,0){\line(1,0){15}}
\multiput(5,0)(5,0){2}{\qbezier[30](0,0)(0,6)(0,12)}
\qbezier(5,0)(7.5,3)(10,0)
\put(9.5,1.1){\vector(3,-1){0}} 
   \put(3.8,2.8){$*$}
   \put(2.1,3.6){\scriptsize $2$}
   \put(3.8,6.7){$*$} 
   \put(2.1,7.5){\scriptsize $1$}
   \put(8.8,4.7){$*$}
   \put(11,5.5){\scriptsize $3$}   
 \end{picture}
\,\Bigr)=
\begin{picture}(10,12)(0,3)
\put(0,0){\line(1,0){10}}
\qbezier[30](5,0)(5,6)(5,12)
\multiput(3.8,-1.1)(0,3.3){4}{$*$}
\put(6,3){\scriptsize $3$}
\put(6,6.3){\scriptsize $2$}
\put(6,9.4){\scriptsize $1$}
\end{picture}
+
\begin{picture}(10,12)(0,3)
\put(0,0){\line(1,0){10}}
\qbezier[30](5,0)(5,6)(5,12)
\multiput(3.8,-1.1)(0,3.3){4}{$*$}
\put(6,3){\scriptsize $2$}
\put(6,6.3){\scriptsize $3$}
\put(6,9.4){\scriptsize $1$}
\end{picture}
+
\begin{picture}(10,12)(0,3)
\put(0,0){\line(1,0){10}}
\qbezier[30](5,0)(5,6)(5,12)
\multiput(3.8,-1.1)(0,3.3){4}{$*$}
\put(6,3){\scriptsize $2$}
\put(6,6.3){\scriptsize $1$}
\put(6,9.4){\scriptsize $3$}
\end{picture}
$$
\vspace{0.2cm}

Obviously, the obtained diagram is always the same for all the shuffles. The calculation of signs shows that this diagram is obtained with the coefficient
${{k+n}\choose{k}}_{(-1)^d}$, where
\begin{gather*}
{ {k+n}\choose{k} }_{1}={ {k+n}\choose{k} };\\
{ {k+n}\choose{k} }_{-1}=
\begin{cases}
0,&\text{$k$ and $n$ are odd;}\\
\left({{\left[\frac{k+n}2\right]}\atop{\left[ \frac k2\right]}}\right), & \text{otherwise}.
\end{cases}
\end{gather*}
This notation is a partial case of the quantum binomial ${{k+n}\choose{k}}_{q}$, when $q$ tends to $\pm 1$. The number ${{k+n}\choose{k}}_{-1}$ is the number
of even shuffles minus the number of odd shuffles. 

The sign of the boundary diagram is obtained by placing $t_{i}$ on the first place of the orienting monomial (this gives a sign) and then removing it.

This complex is a differential bialgebra. One defines the product as a shuffle of active points of diagrams:
$$
\unitlength=0.22em
\begin{picture}(15,12)(0,3)
\put(0,0){\line(1,0){15}}
\multiput(5,0)(5,0){2}{\qbezier[30](0,0)(0,6)(0,12)}
\qbezier(5,0)(7.5,3)(10,0)
\put(9.5,0.9){\vector(3,-1){0}} 
 \end{picture}
\,*\,
\begin{picture}(10,12)(0,3)
\put(0,0){\line(1,0){10}}
\put(5,0){\qbezier[30](0,0)(0,6)(0,12)}
 \put(3.8,5.7){$*$} 
 \end{picture}
=
\begin{picture}(20,12)(0,3)
\put(0,0){\line(1,0){20}}
\multiput(5,0)(5,0){3}{\qbezier[30](0,0)(0,6)(0,12)}
\qbezier(5,0)(7.5,3)(10,0)
\put(9.5,0.9){\vector(3,-1){0}}
  \put(13.8,5.7){$*$} 
 \end{picture}
 +
 \begin{picture}(20,12)(0,3)
\put(0,0){\line(1,0){20}}
\multiput(5,0)(5,0){3}{\qbezier[30](0,0)(0,6)(0,12)}
\qbezier(5,0)(10,5)(15,0)
\put(13.5,1.5){\vector(3,-1){0}}
  \put(8.8,5.7){$*$} 
 \end{picture}
 +
 \begin{picture}(20,12)(0,3)
\put(0,0){\line(1,0){20}}
\multiput(5,0)(5,0){3}{\qbezier[30](0,0)(0,6)(0,12)}
\qbezier(10,0)(12.5,3)(15,0)
\put(14.5,0.9){\vector(2,-1){0}}
  \put(3.8,5.7){$*$} 
 \end{picture}
$$
\vspace{0.08cm}

The coproduct is a coconcatenation:
$$
\unitlength=0.18em
\Delta \Bigl(\,
\begin{picture}(35,12)(0,3)
\put(0,0){\line(1,0){35}}
\multiput(5,0)(5,0){6}{\qbezier[30](0,0)(0,6)(0,12)}
\qbezier(5,0)(10,5)(15,0)
\put(13.5,1.5){\vector(3,-1){0}}
\qbezier(10,0)(15,5)(20,0)
\put(18.5,1.5){\vector(3,-1){0}}

\qbezier(25,0)(27.5,3)(30,0)
\put(29.5,1.1){\vector(3,-1){0}} 
    
   \put(3.8,5){$*$}
 
   \put(3.8,-1.1){$*$} 
 
   \put(23.8,-1.1){$*$}

 \end{picture}
\,\Bigr)
=
\begin{picture}(35,12)(0,3)
\put(0,0){\line(1,0){35}}
\multiput(5,0)(5,0){6}{\qbezier[30](0,0)(0,6)(0,12)}
\qbezier(5,0)(10,5)(15,0)
\put(13.5,1.5){\vector(3,-1){0}}
\qbezier(10,0)(15,5)(20,0)
\put(18.5,1.5){\vector(3,-1){0}}

\qbezier(25,0)(27.5,3)(30,0)
\put(29.5,1.1){\vector(3,-1){0}} 
    
   \put(3.8,5){$*$}
 
   \put(3.8,-1.1){$*$} 
 
   \put(23.8,-1.1){$*$}

 \end{picture}
 {\otimes}
 \begin{picture}(10,12)(0,3)
\put(0,0){\line(1,0){10}}
\end{picture}
+
\begin{picture}(25,12)(0,3)
\put(0,0){\line(1,0){25}}
\multiput(5,0)(5,0){4}{\qbezier[30](0,0)(0,6)(0,12)}
\qbezier(5,0)(10,5)(15,0)
\put(13.5,1.5){\vector(3,-1){0}}
\qbezier(10,0)(15,5)(20,0)
\put(18.5,1.5){\vector(3,-1){0}}
   \put(3.8,5){$*$}
   \put(3.8,-1.1){$*$} 
  \end{picture}
{\otimes}
\begin{picture}(15,12)(0,3)
\put(0,0){\line(1,0){15}}
\multiput(5,0)(5,0){2}{\qbezier[30](0,0)(0,6)(0,12)}
\qbezier(5,0)(7.5,3)(10,0)
\put(9.5,1.1){\vector(3,-1){0}} 
    \put(3.8,-1.1){$*$}
 \end{picture}
 +
  \begin{picture}(10,12)(0,3)
\put(0,0){\line(1,0){10}}
\end{picture}
{\otimes}
\begin{picture}(35,12)(0,3)
\put(0,0){\line(1,0){35}}
\multiput(5,0)(5,0){6}{\qbezier[30](0,0)(0,6)(0,12)}
\qbezier(5,0)(10,5)(15,0)
\put(13.5,1.5){\vector(3,-1){0}}
\qbezier(10,0)(15,5)(20,0)
\put(18.5,1.5){\vector(3,-1){0}}

\qbezier(25,0)(27.5,3)(30,0)
\put(29.5,1.1){\vector(3,-1){0}} 
    
   \put(3.8,5){$*$}
 
   \put(3.8,-1.1){$*$} 
 
   \put(23.8,-1.1){$*$}

 \end{picture}
$$

\section{Derived complexes}\label{s4}
In this section we define complexes $\TstD$, $\TD$, $\TzD$, $\ZZ$ as some quotient-complexes or subcomplexes of $\TdstD$. These complexes were defined before
in~\cite{T5}. 

\begin{definition}\label{d41}
{\rm Complex $\TstD$ is the subcomplex of $\TstD$ spanned by the diagrams without top asterisks
$$
\xymatrix{
\TstD\,\ar@{^{(}->}[r]&\TdstD.
}
$$
the corrisponding diagrams are called} $T_{*}$-diagrams. \nbox
\end{definition}

\begin{definition}\label{d42}
{\rm
Complex $\TD$ is the quotient-complex of $\TdstD$ (and also of $\TstD$) by all diagrams having at least one asterisk.
The diagrams without asterisks are called
} $T$-diagrams. \nbox
\end{definition}

$$
\xymatrix{
\TstD\,\ar@{^{(}->}[rr]\ar@{->>}[rd]&&\TdstD.\ar@{->>}[ld]\\
&\TD&
}
$$

\begin{lemma}\label{l43}
Projection $\TdstD\twoheadrightarrow\TD$ 
 is a quasi-isomorphism. \nbox
\end{lemma}

\noindent {\bf Proof of Lemma~\ref{l43}:} Consider a filtration in $\TdstD$ by the number of active points. The degree zero differential of the spectral sequence
associated to this filtration is $\partial_{v}$. So, the first term (together with the differential on it) of this spectral sequence is the complex $\TD$. By dimensional reasons
all the higher differentials are trivial (for any fixed complexity $i$ the non-trivial groups of the first term are concentrated on only one row). \nbox

\begin{definition}\label{d44}
{\rm
Complex $\TzD$ is a subcomplex of $\TD$ (and also of $\TstD$) spanned by $T$-diagrams without chords joining neighbor active points. The corresponding diagrams 
are called} $T_{0}$-diagrams.\nbox
\end{definition}

\begin{remark}\label{r45}
{\rm It is not evident from the definition that the obtained space forms a subcomplex of $\TD$ or of $\TstD$. But it is so, and moreover $\TzD$ is quasi-isomorphic to $\TstD$,
\cf~\cite{T2,T5}. \nbox
}
\end{remark}

We have a commutative diagram of complexes:
$$
\xymatrix{
\TzD\ar@{^{(}->}[r]^{\sim}\ar@{^{(}->}[d]&\TstD\ar@{^{(}->}[d]\\
\TD&\TdstD\ar@{->>}[l]_{\sim}
}.
\eqno(\numb)\label{eq41}
$$

Recall that differential $\partial$ of complex $\TdstD$ is the sum $\partial=\partial_{h}+\partial_{v}$. We will consider
complex $(\TdstD,\partial_{h})$ whose space is the same space of $T_{*}^*$-diagrams, but the differential is $\partial_{h}$.

\begin{definition}\label{d46}
{\rm Complex $\ZZ$ is a subcomplex of $(\TdstD,\partial_{h})$ spanned by diagrams having only top asterisks
(no chords, no asterisks). \nbox}
\end{definition}

\begin{remark}\label{r47}
{\rm $\ZZ$ is a quotient-complex of $\TdstD$. \nbox}
\end{remark}

\begin{remark}\label{r48}
{\rm All the above complexes are differential bialgebras, product and coproduct being inherited from $\TdstD$.
All the above inclusions and surjections are morphisms of differential bialgebras. \nbox}
\end{remark}

\begin{definition}\label{d49}
{\rm Similarly to Definition~\ref{gtdstd} we define {\it generalized $T_{*}/T/T_{0}$-diagrams} by allowing active points without chords and without asterisks.} \nbox
\end{definition}

\section{Geometric interpretation of the complexes. Main result}\label{s5}
Recall that complex $\TdstD$ is quasi-isomorphic to the $E_{0}$ term of the Vassiliev spectral sequence associated to the space $Emb^+$,
see Proposition~\ref{p31}.

In the same way $\TstD$, resp. $\ZZ$, resp. $(\TdstD,\partial_{h})$ are quasi-isomorphic to the $E_{0}$ terms of the Vassiliev spectrals sequences
associated to the spaces $Emb$, resp. $\Omega Imm\simeq\Omega^2(\R^d\setminus\{ 0\})$, resp. $Emb\times\Omega Imm$. 
For the first and the second complexes, see~\cite{V4}. All these spaces are $H$-spaces what explains why the above complexes are bialgebras. In the table below we resume 
the above correspondence:

\begin{center}
\begin{tabular}{|c|c|c|c|c|c|}
\hline
Space & $Emb^+$ & $Emb$ & $\Omega Imm$ & $Emb\times\Omega Imm$ & 
\parbox{2cm}{{\footnotesize Strata of immersions in} $\Sigma_{d}$}\\
\hline
complex & $\TdstD$ & $\TstD$ & $\ZZ$ & $(\TdstD,\partial_{h})$ & $\TD$ \\
\hline
\parbox{2cm}{\footnotesize complex simplifying computations} & $\TD$ & $\TzD$ & & $\TstD\otimes\ZZ$ & $\TzD\otimes\ZZ$\\
\hline
\end{tabular}
\end{center}

The last entry is not a space, but describes the way how complex $\TD$ was initially obtained.

The main aim of this paper is to prove the following result, which was announced in~\cite[Theorem~14.4]{T5}:

\begin{theorem}\label{t51}
As a differential bialgebra $\TD$ is quasi-isomorphic to $\TzD\otimes\ZZ$ (and to $\TstD\otimes\ZZ$).
\nbox
\end{theorem}

In Section~\ref{s7'}
we give an explicit formula for this quasi-isomorphism.

This result follows from Lemma~\ref{l52} and Theorem~\ref{t53}.

\begin{lemma}\label{l52}
Morphism of differential bialgebras
$$
\mu\colon\TstD\otimes\ZZ\hookrightarrow (\TdstD,\partial_{h})
\eqno(\numb)\label{eq51}
$$
is a quasi-isomorphism. ($\mu$ designates multiplication in the space of $T^*_{*}$-diagrams.) \nbox
\end{lemma}

\noindent {\bf Proof of Lemma~\ref{l52}:} An active point of a $T_{*}^*$-diagram will be called {\it bottom} if it does not contain top astersik, and  is adjacent to some chord or
does contain a bottom asterisk. An active point will be called {\it top} if it contains top asterisks, but does not contain a bottom one and is not adjacent to any chord.
A segment $[t_{i},t_{i+1}]$ is called {\it admissible}, where $t_{i}$ and $t_{i+1}$ are two neighbor active points, if one of the active points is bottom
and another one is top. Consider a filtration in the space of $T_{*}^*$-diagrams by the number of non-admissible segments of neighbor active points.
The first term of the associated spectral sequence (computing the homology of $(\TdstD,\partial_{h})$) together with the first differential is complex $\TstD\otimes\ZZ$.
It is easy to see that other higher differentials are trivial. \nbox

\begin{theorem}\label{t53}
Differential bialgebra $\TdstD$ is isomorphic to $(\TdstD,\partial_{h})$. \nbox
\end{theorem}

In Section~\ref{s7} we will construct this isomorphism:
$$
I\colon (\TdstD,\partial_{h})\stackrel{\simeq}{\longrightarrow}\TdstD.
\eqno(\numb)\label{eq52}
$$

\begin{remark-corollary}\label{c54}{\rm
It was nown that the homology bialgebra of $\TD$ and that of $\TstD^{even}$ are cocommutative for any field of coefficients. The reason for this is that
the duals of $\TD$ and of $\TstD^{even}$ are Hochschild complexes for some operads, and therefore their homology algebras are graded commutative for any ring of coefficients, 
\cf~\cite[Chapitre~IV]{T2}, \cite{T4}, \cite[Part~II]{T5}. On the contrary for $\TstD^{odd}$, the author conjectured that the homology bialgebra of its dual is not
commutative for some ring, \cf~\cite[Conjecture~9.2]{T4}. Theorem~\ref{t51} gives a negative answer to this conjecture. Indeed, as a consequence of this theorem the homology
algebra of the dual of $\TstD$, for any ring of coefficints, is a quotient algebra of the homology algebra of the dual of $\TD$, and hence is commutative. 
One gets as well that the homology bialgebra of $\TstD^{odd}$ is cocommutative for any field of coefficients. \nbox}
\end{remark-corollary}

\section{More algebraic structures}\label{s6}
The aim of this section is to define more algebraic structures on the considered complexes. All the morphismes of complexes respect these structures. Due 
to this it will be sufficient to define isomorphism~\eqref{eq52} on a smaller set of objects.

\begin{definition}\label{d61}
{\rm A {\it divided product} $\langle  D_1,D_2,\dots,D_n\rangle $ of $T_{*}^*/T_{*}$/$T$/$T_0$-diagrams $D_1,D_2,\dots,D_n$ is the sum of those
elements in the shuffle product $D_1*D_2*\ldots *D_n$ that have the left-most point of $D_i$ on the left from the
left-most point of $D_{i+1}$ for all $i=1,\dots,n-1$.} \nbox
\end{definition}

We extend these operations as multilinear operations on the space of $T_{*}$-diagrams (resp. $T$-diagrams and $T_0$-diagrams).

We will denote by $\langle\,\,  \rangle$ the trivial diagram --- the unit of algebras $\TdstD$, $\TstD$, $\TD$ and $\TzD$.

\begin{definition}\label{d63}
Define binary operation
$$
A\vDash B:= (-1)^{|A|+|B|-1}(\partial \langle A,B\rangle -\langle\partial A,B\rangle - (-1)^{|A|}\langle A,\partial B\rangle ).\nboxm
$$
\end{definition}

If $A$ and $B$ are diagrams, then $A\vDash B$ is the sum of all diagrams obtained by gluing the left-most active point of $A$ with the left-most active point of $B$.
Other active points and top asterisks over left-most points are shuffled.

It is easy to see that 
\begin{gather}
A\vDash B=(-1)^{(|A|-1)(|B|-1)}B\vDash A; \label{eq60}\\
(A\vDash B)\vDash C= A\vDash(B\vDash C)\label{eq600}
\end{gather}

\begin{lemma}\label{l64}
\begin{multline*}
\partial\langle A_{1},\ldots, A_{\ell}\rangle = \sum_{i=1}^\ell(-1)^{|A_{1}|+\ldots + |A_{i-1}|} \langle A_{1},\ldots,\partial A_{i},\ldots,A_{\ell}\rangle+\\
+\sum_{i=1}^{\ell-1} (-1)^{|A_{1}|+\ldots + |A_{i}|-1} \langle A_{1},\ldots,A_{i}\vDash A_{i+1},\ldots, A_{\ell}\rangle. \nboxm
\end{multline*}
\end{lemma}

\begin{definition}\label{d65}
{\rm For degree even elements (or for any elements in the case when the characteristics of the main ring is 2)  {\it divided powers} operations are defined as follows}
$$
x^{\langle \ell\rangle}:=\langle \,\underbrace{x,\ldots,x}_{\ell}\,\rangle, \quad \ell=0,1,2,\ldots.\nboxm
$$ 
\end{definition}

\begin{lemma}\label{l65}
Each of the complexes $\TdstD$, $\TstD$, $\TD$, $\TzD$, $(\TdstD,\partial_{h})$, $\ZZ$ is a differential Hopf algebra with 
divided powers. \nbox
\end{lemma}

\noindent {\bf Proof:} It is a direct check of the axioms of Hopf algebra with divided powers, \cf~\cite{Andre}, and also of the axiom
$$
\partial x^{\langle l\rangle} = \partial x \cdot x^{\langle l-1\rangle}.
\eqno(\numb)\label{eq61}
$$
Recall that these axioms say that operations $x\mapsto x^{\langle l\rangle}$ must behave in the same way (with respect to the sum, product, coproduct and differential) as
operations $x\mapsto \frac{x^n}{n!}$ in a differential graded commutative polynomial bialgebra defined over $\Q$. For instance, let us prove~\eqref{eq61}. Note, that $x$ 
is of even degree (or of any degree if characteristics is 2). It can be easily seen that $x\vDash x=0$ in the above case. Applying Lemma~\ref{l64} one gets the result. \nbox

\begin{definition}\label{d66}
{\rm
Let $A_{1},\ldots,A_{\ell}$ be generalized $T_{*}^*$-diagrams, and $D$ be a generalized $T_{*}$-diagram with exactly $\ell$ 
points.  Define
$$
\langle A_{1},\ldots, A_{\ell} | D\rangle
\eqno(\numb)\label{eq62}
$$
as follows. Suppose $t_{1}<t_{2}<\ldots<t_{\ell}$ are the active points of $D$.
If for some $i=1,\ldots,\ell$ the $i$-th active point $t_{i}$ of $D$  and the left-most point of $A_{i}$  contain both  bottom asterisk, then~\eqref{eq62} is defined as zero.
Otherwise~\eqref{eq62} is the sum of the diagrams obtained by the following recipe

1) chose a digram $\tilde A$ in the sum $\langle A_{1},\ldots,A_{\ell}\rangle$ (diagram $\tilde A$ is defined by a shuffle of active points of the $A_{i}$, 
$i=1,\ldots,\ell$);

2) glue  (in the prescribed order) the active points $t_{1},\ldots,t_{\ell}$ of $D$ to the left-most active points of $A_{1},\ldots,A_{\ell}$ in $\tilde A$. 

3) orientation of the obtained diagram is given as follows: Orient $D$ in the way that on  the first $\ell$ places its orienting monomial has $t_{1},\ldots,t_{\ell}$
(this gives a sign). Then we remove this part from the orienting monomial of $D$, and concatenate the orienting monomial of $\tilde A$ and the rest of the orienting
monomial of $D$. \nbox
}
\end{definition}

Recall, that 
$$
\partial D=\partial_{h}D=\sum_{i=1}^{l-1}\partial_{i}D,
$$
where $\partial_{i}$ is the gluing of $t_{i}$ with $t_{i+1}$.

The following two lemmas are easy to verify.

\begin{lemma}\label{l67}
\begin{multline*}
\partial\langle A_{1},\ldots, A_{\ell}|D\rangle = \sum_{i=1}^\ell(-1)^{|A_{1}|+\ldots + |A_{i-1}|} \langle A_{1},\ldots,\partial A_{i},\ldots,A_{\ell}
|D\rangle+\\
+\sum_{i=1}^{\ell-1} (-1)^{|A_{1}|+\ldots + |A_{i}|-1} \langle A_{1},\ldots,A_{i}\vDash A_{i+1},\ldots, A_{\ell}|\partial_{i}D\rangle. \nboxm
\end{multline*}
\end{lemma}

\begin{lemma}\label{l68}
Assertions of Lemmas~\ref{l64} and~\ref{l67} hold if $\partial$ is replaced by $\partial_{h}$. \nbox
\end{lemma}

\section{Isomorphism $I\colon (\TdstD,\partial_{h})\to\TdstD$}\label{s7}
Now we are ready to describe isomorphism~\eqref{eq52}. Define $I$ to be the identity on $\TstD\subset (\TdstD,\partial_{h})$.

Let us define $I$ on $\ZZ\subset(\TdstD,\partial_{h})$. Denote by $Z_{k}$ the diagram:
$$
Z_{k}=
\unitlength=0.25em
\begin{picture}(10,15)(0,3)
\put(0,0){\line(1,0){10}}
\qbezier[30](5,0)(5,7.5)(5,15)

\put(4,2.2){$*$}
\put(4,8.8){$*$}
\put(4,12.1){$*$}
\put(6,2.5){\scriptsize $k$}
\put(6.4,5){\scriptsize $\vdots$}
\put(6,9.4){\scriptsize $2$}
\put(6,12.7){\scriptsize $1$}
\end{picture}
\eqno(\numb)\label{eq71}
$$
\vspace{0.1cm}

\noindent oriented by the monomial
$$
t_{1}(\alpha^*{}^1_1a_1^1)(\alpha^*{}^2_1a_1^2)\ldots(\alpha^*{}^k_1a_1^k).
$$
Obviously, $Z_{k}$ is a cycle in $\ZZ$:
$$
\partial_{h}Z_{k}=0.
\eqno(\numb)\label{eq72}
$$
Note also, that
$$
\partial Z_{k}=\partial_{v}Z_{k}=(-1)^{dk}Z_{k-1}\vDash {\bigstar},
\eqno(\numb)\label{eq73}
$$
where $\bigstar$ designates the diagram 
$$\bigstar=\,
\unitlength=0.2em
\begin{picture}(10,7)
\put(0,0){\line(1,0){10}}
\qbezier[12](5,0)(5,3.5)(5,7)
\put(3.82,-1.3){$*$}
\end{picture}
$$
 oriented by $t_{1}\alpha_{1}^*$.

Denote by $\hat Z_{k}$ the diagram
$$
\hat Z_{k}=
\unitlength=0.3em
\begin{picture}(35,13)(0,0)
  \put(0,0){\line(1,0){35}}
  \qbezier(5,0)(7.5,4)(10,0)
  \put(10,0){\vector(2,-3){0}}
  \qbezier(5,0)(10,7)(15,0)
  \put(15,0){\vector(2,-3){0}}
  \qbezier(5,0)(15,12)(25,0)
  \put(25,0){\vector(2,-3){0}}
  \qbezier(5,0)(17.5,15)(30,0)
  \put(30,0){\vector(2,-3){0}}
  \put(16,1.2){\Large $\ldots$}
  \put(4.1,-2.7){$0$}
  \put(9.1,-2.7){$1$}
  \put(14.1,-2.7){$2$}
   \put(29.1,-2.7){$k$}
\end{picture}
\eqno(\numb)\label{eq74}
$$
\vspace{0.1cm}

\noindent oriented by the monomial $t_{0}(\alpha_{01}t_1)\ldots(\alpha_{0k}t_k)$.

It is easy to see that
$$
\partial\hat Z_{k}=(-1)^{d-1}\bigstar\vDash Z_{k-1}.
\eqno(\numb)\label{eq75}
$$

\medskip

Define
$$
I(Z_{k}):=\sum_{i=0}^kZ_{k-i}\vDash\hat Z_{i}.
\eqno(\numb)\label{eq76}
$$

$$
I(Z_{2})=
\unitlength=0.2em
I \Bigl(\,
\begin{picture}(10,10)(0,3)
\put(0,0){\line(1,0){10}}
\put(5,0){\qbezier[15](0,0)(0,5)(0,10)}
      \put(3.8,2.5){$*$}
   \put(3.8,6){$*$} 
    \end{picture}
\,\Bigr)
=
\begin{picture}(10,10)(0,3)
\put(0,0){\line(1,0){10}}
\put(5,0){\qbezier[15](0,0)(0,5)(0,10)}
      \put(3.8,2.5){$*$}
   \put(3.8,6){$*$} 
\end{picture}
+
\begin{picture}(15,10)(0,3)
\put(0,0){\line(1,0){15}}
\multiput(5,0)(5,0){2}{\qbezier[15](0,0)(0,5)(0,10)}
\put(10,0){\vector(1,-1){0}}
      \put(3.8,3.5){$*$}
 \qbezier(5,0)(7.5,3)(10,0)
\end{picture}
+
\begin{picture}(20,10)(0,3)
\put(0,0){\line(1,0){20}}
\multiput(5,0)(5,0){3}{\qbezier[15](0,0)(0,5)(0,10)}
 \qbezier(5,0)(7.5,3)(10,0)
 \put(10,0){\vector(1,-1){0}}
\qbezier(5,0)(10,5)(15,0)
\put(15,0){\vector(1,-1){0}}
\end{picture}
$$
\vspace{0.1cm}

$Z_{0}=\hat Z_{0}$ designates in~\eqref{eq76} and in the sequel the generalized diagram
$$
Z_{0}=\hat Z_{0}=
\unitlength=0.2em
\begin{picture}(10,8)(0,2)
\put(0,0){\line(1,0){10}}
\put(5,0){\qbezier[12](0,0)(0,4)(0,8)}
\end{picture}
$$
with the only active point, without asterisks and without chords.

\begin{lemma}\label{l71}
$\partial I(Z_{k})=0=I(\partial_{h}Z_{k}). \nboxm$
\end{lemma}

\noindent {\bf Proof:} Consequence of~\eqref{eq73}, \eqref{eq75} and also of proprieties~\eqref{eq60}, \eqref{eq600} of $\vDash$. \nbox

\bigskip
Any diagram of $\ZZ$ is of the form $\langle Z_{k_{1}},\ldots,Z_{k_{\ell}}\rangle$. Define 
$$
I\bigl|_{\ZZ}\colon \ZZ\hookrightarrow\TdstD
\eqno(\numb)\label{eq77}
$$
by
$$
I\colon\langle Z_{k_{1}},\ldots,Z_{k_{\ell}}\rangle\mapsto \langle IZ_{k_{1}},\ldots,IZ_{k_{\ell}}\rangle
\eqno(\numb)\label{eq78}
$$

\begin{lemma}\label{l72}
\eqref{eq77}  is a morphism of differential Hopf algebras with divided powers. \nbox
\end{lemma}

\noindent {\bf Proof:} The most difficult  is to prove that it is a morphisme of complexes. As a consequence of Lemmas~\ref{l64}, \ref{l68}, \ref{l71},
it is sufficient to prove that
$$
I(Z_{a}\vDash Z_{b})=I(Z_{a})\vDash I(Z_{b}),
$$
for any integers $a,b>0$.

Note, that $Z_{a}\vDash Z_{b}={{a+b}\choose a}_{(-1)^d}Z_{a+b}$, and $\hat Z_{a}\vDash \hat Z_{b}={{a+b}\choose a}_{(-1)^d}\hat Z_{a+b}$.
The rest is a calculation. \nbox

\bigskip

Now, we are ready to define $I$ on the whole complex $(\TdstD,\partial_{h})$. Let $D$ be any $T_{*}^*$-diagram with $\ell$ active points, then
$D$ can be uniquely presented in the form
$$
D=\langle Z_{k_{1}},\ldots,Z_{k_{\ell}}|D_{bottom}\rangle,
$$
where $k_{i}\geq 0$, $i=1,\ldots,\ell$, and $D_{bottom}$ is a generalized $T_{*}$-diagram, see Definition~\ref{d49}. In fact $D_{bottom}$ is obtained
from $D$ by removing all top asterisks; $k_{i}$ is  the number of top asterisks over $i$-th active point of $D$.

We define
$$
I(D)=I\langle Z_{k_{1}},\ldots,Z_{k_{\ell}}|D_{bottom}\rangle:=\langle IZ_{k_{1}},\ldots,IZ_{k_{\ell}}|D_{bottom}\rangle.
$$

Note for example, that $I(D)=D$, if $D\in\TstD$.

\begin{proposition}\label{p73}
The above map
$$
I\colon (\TdstD,\partial_{h})\to \TdstD
$$
is an isomorphism of differential Hopf algebras with divided powers. Moreover $I$ respects operation $\vDash$. \nbox
\end{proposition}

\noindent {\bf Proof:} It is easy to find a base in the space $\TdstD$ in which the map $I$ is triangular with units on the main diagonal.
So, $I$ is an isomorphism of linear spaces. Lemmas~\ref{l67}, \ref{l68}, \ref{l71} imply that $I$ is a morphism of complexes. The rest is a tedious
check. \nbox

\bigskip

\section{\bf Proof of Theorem~\ref{t51}}\label{s7'} 
Let $\hat\ZZ$ designates a subcomplex of $\TD$ spanned by elements $\langle\hat Z_{k_{1}},\ldots,\hat Z_{k_{\ell}}\rangle$, $k_{i}>0$, $i=1\ldots \ell$.
Obviously,
$$
\hat I\colon Z\to\hat Z
$$
is an isomorphism, where
$$
\hat I\colon \langle Z_{k_{1}},\ldots,Z_{k_{\ell}}\rangle\mapsto \langle \hat Z_{k_{1}},\ldots,\hat Z_{k_{\ell}}\rangle.
$$
By abuse of the language let $\hat I$ designate the composition map $\hat I\colon \ZZ\otimes\TzD\to\TD$:
$$
\xymatrix{
\ZZ\otimes\TzD\ar[r]^{\hat I\otimes \id}&\hat\ZZ\otimes\TzD\ar[r]^-\mu&\TD,
}
$$
where $\mu$ is the shuffle product in $\TD$.

$\hat I$ is the desired morphism. To see that $\hat I$ is a quasi-isomorphism we decompose it in a composition of quasi-isomorphisms:
$$
\xymatrix{
\ZZ\otimes\TzD\,\ar@{^{(}->}[r]^\sim&\ZZ\otimes\TstD\,\ar@{^{(}->}[r]^-{\mu}&(\TdstD,\partial_{h})\ar[r]^-{I}&\TdstD\ar@{->>}[r]^\sim&\TD.
}
$$
The first map is quasi-isomorphism by Remark~\ref{r45}; the second one --- by Lemma~\ref{l52}; the third one --- by Proposition~\ref{p73};
the forth one --- by Lemma~\ref{l43}. \nbox 

\section{Inverse map $I^{-1}$}\label{s8}
The inverse map
$
I^{-1}\colon\TdstD\to (\TdstD,\partial_{h})
$
can be described explicitly. 

\begin{proposition}\label{p81}
The inverse map $I^{-1}$ is described as follows.

1) $I^{-1}(Z_{k})=\sum_{i=0}^k(-1)^{i+d\frac{i(i-1)}2}Z_{k-i}\vDash\hat Z_{i}$, $k=0,1,2,\ldots$.

2) For any diagram $D\in\TdstD$, 
$I^{-1}(D)=I^{-1}\langle Z_{k_{1}},\ldots,Z_{k_{\ell}}|D_{bottom}\rangle=\langle I^{-1}Z_{k_{1}},\ldots,I^{-1}Z_{k_{\ell}}|D_{bottom}\rangle$,
$k_{i}\geq 0$, $i=1,\ldots,\ell$, where $D_{bottom}$ and $k_{i}$'s are determined by $D$ as in Section~\ref{s7}. \nbox
\end{proposition}

\section{Upper diagonal of the  spectral sequences}\label{s9}
For a given complexity $i$, the number $j$ of geometrically distinct points of a $T$-diagram can be in the range:
$$
i+1\leq j\leq 2i.
$$
As a consequence of~\eqref{eq32}, Lemma~\ref{l43} and Remark~\ref{r45}, the first term of the Vassiliev spectral sequence associated to 
$Emb^+$ (and to $Emb$) is concentrated in the second quadrant
\begin{figure}[!ht]
\begin{center}
\unitlength=0.5em
\begin{picture}(20,30)
\put(0,10){\vector(1,0){20}}
\put(10,0){\vector(0,1){30}}
\put(0,30){\line(1,-2){10}}
\put(5,30){\line(1,-4){5}}
\put(19,8.5){$p$}
\put(10.5,29){$q$}
\put(8,14){\line(1,0){1}}
\put(6,18){\line(1,0){2}}
\put(4,22){\line(1,0){3}}
\put(2,26){\line(1,0){4}}
\end{picture}
\end{center}
\caption{First term of the Vassiliev spectral sequence}\label{SpSeq}
\end{figure}
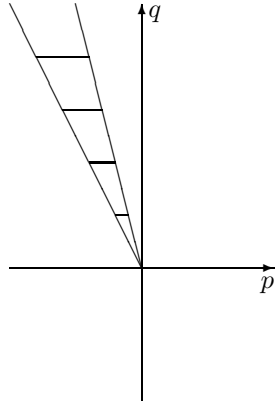
between two half-lines:
$$
\begin{tabular}{ll}
$q=-(d-2)p$&\text{lower line;}\\
$q=-(d-1)p-1$&\text{upper line.}
\end{tabular}\eqno(\numb\label{eq91})
$$
The groups of the lower half-line $(j=2i)$ form the bialgebra $\BB$ (resp. $\BB_{0}$) of chord diagrams, for odd $d$, and some its 
non-trivial super-analogue $\tilde\BB$ (resp. $\tilde\BB_{0}$), for even $d$. The upper-diagonal groups are described by the following theorem.

\begin{theorem}\label{t91}{\rm \cite{T5}}
{\rm (i)} Homology groups of $\TzD^d$ are trivial in the bigradings $j=i+1$.

{\rm (ii)} Homology groups of $\TD^d$ in the bigradings $(i,i+1)$ are as follows
\begin{gather*}
\text{$d$ is even}\quad\begin{cases}
\Z,&i=1;\\
\Z_p,&\text{$i=p^k$, where $p$ is any prime, $k\in\N$;}\\
0,&\text{otherwise.}
\end{cases}
\\
\text{$d$ is odd}\quad
\begin{cases}
\Z,&i=1,\,\, 2;\\
\Z_p,&\text{$i=2p^k$, where $p$ is any prime, $k\in\N$;}\\
0,&\text{otherwise.}
\end{cases}
\end{gather*}
As a generator one can choose diagram $\hat Z_{i}$, see~\eqref{eq74}. \nbox
\end{theorem}

The proof of (i) is very simple, but the proof of (ii) was more tedious in~\cite{T5}. We used there the Dyer-Lashof-Cohen operations on the dual complex 
(which is Hochschild complex of $(d-1)$-Poisson algebras operad) to describe explicitly the dual cycles. Theorem~\ref{t51} together with (i) provides another proof
of (ii).

\begin{theorem}\label{t93}
All the non-trivial upper-diagonal groups $E_{1}^{p,-(d-1)p-1}$ of the Vassiliev spectral sequence computing the cohomology of $Emb^+_{d}$ survive the higher differentials:
$E_{\infty}^{p,-(d-1)p-1}=E_{1}^{p,-(d-1)p-1}$, and define non-trivial cohomology classes of the space $Emb_{d}^+$, $d\geq 3$. \nbox
\end{theorem}

Recall that $Emb^+_{d}\simeq Emb\times\Omega^2 S^{d-1}$. The homology of $\ZZ$ is exactly $H^*(\Omega^2 S^{d-1})$, see~\cite{V4}\footnote{See also~\cite{F,Vain,V0,Ma,T5}
where the homology groups of $\ZZ$ are computed.}. So, Theorem~\ref{t51}
says that for any field of coefficients $E_1^{*,*}(Emb_{d}^+)\simeq E_1^{*,*}(Emb_{d})\otimes H^*(\Omega^2 S^{d-1})$. It seems natural
that this splitting would take place for higher terms $E_{r}$ of the Vassiliev spectral sequence: 
$E_{r}^{*,*}(Emb_{d}^+)\simeq E_r^{*,*}(Emb_{d})\otimes H^*(\Omega^2 S^{d-1})$. But we have no natural morphism between 
spectral sequences, so we can not prove this splitting for higher~$r$.

Our proof of Theorem~\ref{t93} is an explicit geometric calculation of the higher differentials. The technique that we use is that of~\cite{V5,V6}.
The main idea is as follows. One considers the geometrical border of a higher differential. This defines a cycle in the smaller filtration term of the resolved discriminant.
It turns out that this cycle is always trivial. To see this one finds a proper chain in the smaller filtration term, whose border is exactly this cycle. Doing it step by step
one gets a cycle in the resolved discriminant. Once this cycle is obtained, one considers its projection (from the resolved discriminant) to the discriminant itself.
The linking number with the cycle (which is the image of the projection) is exactly the desired cohomology class. To be able to compute the linking number, one finds a 
chain in the space $Emb^+_{d}$ of non-singular maps whose border is the latter cycle. Now, the cohomology class is defined as the intersecting number with the obtained 
chain.

We will not describe all the steps of these calculations in the case of the upper diagonal cocycles, but present straight away the obtained infinite-dimensional chains.
So, we leave it as an exercise to the reader to see that the technique of V.~Vassiliev in~\cite{V5,V6} provides exactly the result described below.

\medskip

Before describing this result let us consider a simpler situation --- that of the space $\Omega Imm$. The first term of the Vassiliev spectral sequence
associated to this space is concentrated in the same sector, see Figure~\ref{SpSeq}, and is formed by the homology groups of the complex $\ZZ$.
The lower diagonal groups correspond to bigradings $j=2i$, \ie to diagrams with one asterisk over each active point:
$$
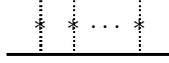
\unitlength=0.25em
\begin{picture}(25,8)
\put(0,0){\line(1,0){25}}
\put(5,0){\qbezier[12](0,0)(0,4)(0,8)}
\put(10,0){\qbezier[12](0,0)(0,4)(0,8)}
\put(20,0){\qbezier[12](0,0)(0,4)(0,8)}
 \put(4,3.5){$*$}
 \put(9,3.5){$*$}
 \put(12.5,4){$\ldots$} 
 \put(19,3.5){$*$}
\end{picture}
$$
The upper diagonal (j=i+1) corresponds to diagrams of $\ZZ$ with only one active point:
$$
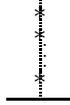
\unitlength=0.25em
\begin{picture}(10,18)(0,0)
\put(0,0){\line(1,0){10}}
\qbezier[30](5,0)(5,7.5)(5,15)

\put(4,2.2){$*$}
\put(4,8.8){$*$}
\put(4,12.1){$*$}
\put(5.3,4.8){\scriptsize $\vdots$}
\end{picture}
$$
The homology groups of $\ZZ$ in bigradings $j=i+1$ are described by (ii) of Theorem~\ref{t91}. So, any
diagram $Z_{i}$ defines a cycle (over $\Z$) in the filtration $\sigma_{i}\setminus\sigma_{i-1}$ which is untrivial only if 
$i=p^k$ ($d$ even) or $i=2p^k$ ($d$ odd). Denote by $V(Z_{i})$ the cocycle of $\Omega Imm$ obtained from
$Z_{i}$ by the described above procedure.

\begin{theorem}\label{t94}
The value of the cocycle $V(Z_{i})$ on a generic cycle $C$ in $\Omega Imm_{d}$ of dimension $i(d-2)-1$ is equal to the number of points $f\in\Omega Imm_{d}$ of $C$
(counted with appropriate signs), such that for some $t_{0}\in\R^1$, and some $1>a_{1}>a_{2}>\ldots>a_{i}>0$ one has
$$
f'_{t}(a_s,t_0)=(-1)^{s-1}\lambda_{s}\partial/\partial x_{1},
$$
for some $\lambda_{s}>0$, $s=1\ldots i$. \nbox
\end{theorem}

The corresponding infinite dimensional chain (with which one takes intersecting number, see above) can be represented by the diagram:

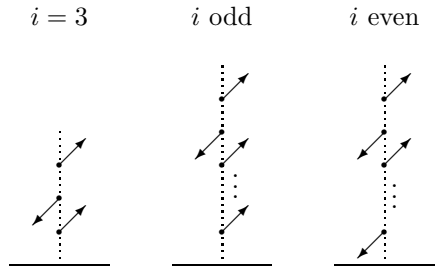
\unitlength=0.25em
\begin{center}
\begin{tabular}{ccccc}
$i=3$&\quad&$i$ odd&\quad&$i$ even\\
\\
\begin{picture}(15,20)
\put(0,0){\line(1,0){15}}
\qbezier[20](7.5,0)(7.5,10)(7.5,20)
\put(7.5,5){\circle*{0.7}}
\put(7.5,10){\circle*{0.7}}
\put(7.5,15){\circle*{0.7}}
\put(7.5,5){\vector(1,1){4}}
\put(7.5,10){\vector(-1,-1){4}}
\put(7.5,15){\vector(1,1){4}}
\end{picture}
&&
\begin{picture}(15,30)
\put(0,0){\line(1,0){15}}
\qbezier[30](7.5,0)(7.5,15)(7.5,30)
\put(7.5,5){\circle*{0.7}}
\put(7.5,15){\circle*{0.7}}
\put(7.5,20){\circle*{0.7}}
\put(7.5,25){\circle*{0.7}}
\put(7.5,5){\vector(1,1){4}}
\put(9,10){$\vdots$}
\put(7.5,15){\vector(1,1){4}}
\put(7.5,20){\vector(-1,-1){4}}
\put(7.5,25){\vector(1,1){4}}
\end{picture}
&&
\begin{picture}(15,30)
\put(0,0){\line(1,0){15}}
\qbezier[30](7.5,0)(7.5,15)(7.5,30)
\put(7.5,5){\circle*{0.7}}
\put(7.5,15){\circle*{0.7}}
\put(7.5,20){\circle*{0.7}}
\put(7.5,25){\circle*{0.7}}
\put(7.5,5){\vector(-1,-1){4}}
\put(8.5,8.5){$\vdots$}
\put(7.5,15){\vector(1,1){4}}
\put(7.5,20){\vector(-1,-1){4}}
\put(7.5,25){\vector(1,1){4}}
\end{picture}
\end{tabular}
\end{center}

where the arrow 
\begin{picture}(4,4)
\put(0,0){\vector(1,1){4}}
\end{picture}
 (resp. 
 \begin{picture}(4,4)
\put(4,4){\vector(-1,-1){4}}
\end{picture}\,)
 means that the derivative vector $f'_{t}$ is directed \lq\lq up" --- in the direction of the axis $x_{1}$ (resp. \lq\lq down" --- in the opposite 
direcion to this axis). We suppose that the direction $\partial/\partial x_{1}$ does not coincide with the direction of the fixed linear map.

The directions \lq\lq up" and \lq\lq down" alternate. It is done delibirately in order the corresponding infinite dimensional chain would have a border only in the 
complement of $\Omega Imm_{d}$.

\bigskip

Let us come back to the space $Emb^+$. Denote by $V(I(Z_{i}))$, $i=p^k$ ($d$ even), $i=2p^k$ ($d$ odd), cocycles of $Emb^+$ obtained by the similar procedure, \cf~\cite{V5,V6}.

\begin{theorem}\label{t95}
The value of the cocycle $V(I(Z_{i}))$ on a generic cycle $C$ in $Emb^+_{d}$ of dimension $i(d-2)-1$ is equal to the number of points $f\in Emb^+_{d}$ of $C$
(counted with appropriate signs) such that for some integer $\ell$, $0\leq \ell\leq i$, there exist points $t_0<t_1<\ldots<t_{i-\ell}$ in $\R$, and points $a_1>a_2>\ldots>a_{\ell}$
in $(0,1)$ for which the following holds:
\begin{gather*}
f'_t(a_s,t_0)=(-1)^{s-1}\lambda_{s}\partial/\partial x_{1}, \quad 1\leq s\leq \ell,  \\
f(0,t_{r})-f(0,t_{0})=(-1)^{r+\ell-1}\mu_{r}\partial/\partial x_{1}, \quad 1\leq r\leq i-\ell,
\end{gather*}
for some $\lambda_{s}>0$, $1\leq s\leq \ell$; $\mu_{r}>0$, $1\leq r\leq i-\ell$. \nbox
\end{theorem}

The corresponding infinite dimensional chain in $Emb^+$ (with wich one counts the intersecting number) can be represented as a sum of $(i+1)$ diagrams. For example, for $i=3$
this sum is:
$$
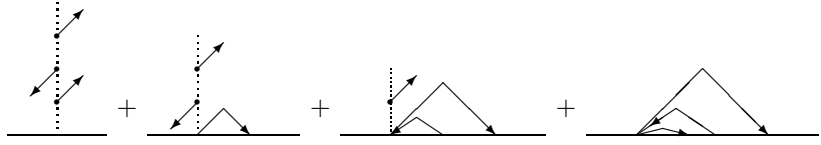
\unitlength=0.25em
\begin{picture}(15,20)(0,3)
\put(0,0){\line(1,0){15}}
\qbezier[20](7.5,0)(7.5,10)(7.5,20)
\put(7.5,5){\circle*{0.7}}
\put(7.5,10){\circle*{0.7}}
\put(7.5,15){\circle*{0.7}}
\put(7.5,5){\vector(1,1){4}}
\put(7.5,10){\vector(-1,-1){4}}
\put(7.5,15){\vector(1,1){4}}
\end{picture}
\,
+
\,
\begin{picture}(23,15)(0,3)
\put(0,0){\line(1,0){23}}
\qbezier[15](7.5,0)(7.5,7.5)(7.5,15)
\put(7.5,5){\circle*{0.7}}
\put(7.5,10){\circle*{0.7}}
\put(7.5,5){\vector(-1,-1){4}}
\put(7.5,10){\vector(1,1){4}}
\put(7.5,0){\line(1,1){4}}
\put(11.5,4){\vector(1,-1){4}}
\end{picture}
\,
+
\,
\begin{picture}(31,10)(0,3)
\put(0,0){\line(1,0){31}}
\qbezier[15](7.5,0)(7.5,5)(7.5,10)
\put(7.5,5){\circle*{0.7}}
\put(7.5,5){\vector(1,1){4}}
\put(7.5,0){\line(1,1){8}}
\put(15.5,8){\vector(1,-1){8}}
\put(15.5,0){\line(-3,2){4}}
\put(11.5,2.66){\vector(-3,-2){4}}
\end{picture}
\,+\,
\begin{picture}(35,10)(0,3)
\put(0,0){\line(1,0){35}}
\put(7.5,0){\line(1,1){10}}
\put(17.5,10){\vector(1,-1){10}}

\put(7.5,0){\line(3,2){6}}
\put(13.5,4){\line(3,-2){6}}
\put(9.5,1.33){\vector(-3,-2){0}}

\put(11.5,1){\vector(4,-1){4}}
\put(11.5,1){\line(-4,-1){4}}
\end{picture}
$$
\vspace{0.2cm}

A broken arrow from $t_{i}$ to $t_{j}$ means that the point $f(t_{j})$ is \lq\lq over" the point $f(t_{i})$.

Note that the first summand of this chain is the chain from Theorem~\ref{t94}. So, the restriction of the cocycle $V(I(Z_{i}))$ on the space $\Omega Imm_{d}\subset Emb^+_{d}$
is exactly $V(Z_{i})$.

\medskip

The above chains are defined over $\Z$, but if $i\neq 1$ (any $d$) and $i\neq 2$ (odd $d$), then the intersection number with any $\Z$-cycle is always zero. One should consider cycles
over the corresponding cyclic group $\Z_{p}$ to obtain non-trivial intersections.

\section{Approach of T.~Goodwillie. Sinha's spectral sequence}\label{s10}
There is another and absolutely different approach to studying the spaces of embeddings. This approach is an application of the \lq\lq Calculus of Functors" --- theory
developped by T.~Goodwillie. Briefly speaking in this approach one \lq\lq approximates" the space of knots by homotopy limits of digrams of maps.

Dev Sinha used this method and showed that the space $Emb_{d}$, $d\geq 4$, is the homotopy totalization of the cosimplicial space $[C'_{d}]^\bullet$:
$$
\xymatrix{
 [C'_{d}]^0\ar@<0.7ex>[r]\ar@<-0.7ex>[r] & [C'_{d}]^1\ar[l]\ar@<1ex>[r]\ar@<-1ex>[r]\ar[r] &  [C'_{d}]^2\ar@<0.5ex>[l]\ar@<-0.5ex>[l]
\ar@<0.4ex>[r]\ar@<1.2ex>[r]\ar@<-0.4ex>[r]\ar@<-1.2ex>[r]&\dots
},
$$
whose
$n$-th component $[C'_{d}]^n$ is the configuration space of distinct points with  unit tangent vectors in these points:
$$
[C'_{d}]^n=\left\{ \, (x_0,\ldots,x_{n+1};v_0,\ldots,v_{n+1})\, \Biggl| \, \text{\parbox{6.5cm}{$x_{i}\in [0,1]\times \R^{d-1}$, $v_{i}\in S^{d-1}$, $i=0\ldots n{+}1$; 
$x_{0}{=}(0,\bar 0)$, $x_{n+1}{=}(1,\bar 0)$, $v_{0}{=}v_{n+1}{=}(1,\bar 0)$; $x_{i}\neq x_{j}$, $0\leq i\neq j\leq n+1$}} \, \right\}.
$$
To be precise $[C'_d]^n$ is some compactification of the above space, \cf~\cite{Sinha1}.

Coface maps
$$
d_{i}\colon [C'_{d}]^k\to [C'_{d}]^{k+1}, \quad i=0\ldots k{+}1,
$$ 
are defined as doubling of the $i$-th point in the configuration in the direction of the $i$-th tangent vector $v_{i}$.

Codegeneracy maps
$$
s_{i}\colon [C'_{d}]^k\to [C'_{d}]^{k-1}, \quad i=1\ldots k,
$$
are defined as forgetting of the $i$-th point.

\medskip

In~\cite[Section~7]{Sinha1}, D.~Sinha defined a spectral sequence computing the cohomology groups of the homotopy totalization of  $[C'_{d}]^\bullet$ (\ie of the space 
$Emb_{d}$). The first term $E_{1}$ of  Sinha's spectral sequence is exactly the complex $\TstD^d$. Sinha's $E_{1}$ term is the normalized part of the simplicial
algebra $H^*([C'_{d}]^\bullet)$ formed by the cohomology algebras  
$H^*([C'_{d}]^k), \, k=0,1,2,\ldots$ Differential in $E_{1}$ is the alternated sum of $d_{i}^*$ (cohomology maps
induced by cofaces~$d_{i}$).  "Normalized part" is the quotient complex (it is always isomorphic to the initial one) of the initial complex formed by
$\oplus_{n\geq 0}H^*([C'_{d}]^{n})$. The quotient is taken by the sum of images of maps $s_{i}^*$.

In~\cite{Sinha2} D.~Sinha considers a similar cosimplicial space $[C_{d}]^\bullet$, which is a subspace of  $[C'_{d}]^\bullet$, setting all the unit tangent vectors
$v_{i}$ to be $(1,\bar 0)$. Using homotopy limits' proprieties, he shows that the homotopy totalization of  $[C_{d}]^\bullet$ is the homotopy fiber of the inclusion
$Emb_{d}\hookrightarrow Imm_{d}$, \ie our space $Emb^+_{d}$. On the other hand, the first term of the analogous spectral sequence computing the cohomology
of the totalization is the complex $\TD$, \cf~\cite[Section~7]{Sinha2}.

\smallskip

One has the following result.

\begin{theorem}\label{t101}
{\rm P.~Lambrechts, I.~Volic~\cite{LV}} For $d\geq 4$, D.~Sinha's spectral sequences computing the cohomology groups of $Emb_{d}$ and of
$Emb_{d}^+$ collapse over $\Q$ in the second term. \nbox
\end{theorem}

\begin{corollary}\label{c102}
For $d\geq 4$, the Vassiliev spectral sequences computing the (co)homology of $Emb_{d}$ and $Emb_{d}^+$ collapse over $\Q$ in the first term. \nbox
\end{corollary}

The direct sum of all the components of a simplicial module form a complex, where the differential is the alternated sum of face maps.
The homology of  this complex is usually designated by $\pi_{*}$ and is called the {\it homotopy} of the simplicial module, 
\cf~\cite{Dwyer, Fresse1,Fresse2}.

\smallskip

The direct sum of all the components of a simplicial algebra forms a differential graded algebra. Multiplication being defined via the map analogous to the Eilenberg-Zilber morphism. 
Recall that the Eilenberg-Zilber map is the natural quasi-isomorphism from the tensor product of the singular chains complex  of a space $X$ with singular chains complex  of a space $Y$
$$
\xymatrix{
EZ\colon S_{*}X\otimes S_{*}Y\,\ar@{^{(}->}[r]^-\sim& S_{*}(X\times Y)}
\eqno(\numb)\label{eq101}
$$
to the singular chains complex of the product $X\times Y$,\cf~\cite[Chapter~8]{Mac}, \cite[Section~29]{May1}. 
Map $EZ$ is also called shuffle map. 
Geometrically it corresponds to natural simplicial subdivisions
of prisms $\Delta^k\times \Delta^\ell$ which are products of two simplices $\Delta^k$ and $\Delta^\ell$. 

Analogously one defines a quasi-isomorphism from the tensor product $C_{*}(V)\otimes C_{*}(W)$ of complexes formed by two simplicial modules $V$ and $W$
to the complex formed by the simplicial module which is the component-wise
tensor product of $V$ and $W$: 
$$
(V\otimes W)_{n}=V_{n}\times W_{n}.
$$

Multiplication in a simplicial algebra $A$ is defined as the composition of the map analogous to~\eqref{eq101} (which sends tensor square of the complex formed by a simplicial algebra to the
component-wise tensor square), and component-wise multiplication $\mu$, \cf~\cite{Dwyer, Fresse1,Fresse2}:
$$
\xymatrix{
C_{*}(A)\otimes C_{*}(A)\,\ar@{^{(}->}[r]^-{EZ}& C_{*}(A\otimes A)\ar[r]^-\mu& C_{*}(A).}
$$
This product is associative and graded commutative (if $A$ is component-wise associative and commutative) and  induces an associative and graded commutative 
product on the normalized part $N_{*}(A)$ of simplicial algebra $A$.

The following assertion is easy to verify.

\begin{lemma}\label{l101}
Multiplication in $\TD$, $\TstD$ which was defined in Section~\ref{s3} is exactly the multiplication induced by the simplicial algebra structure on $H^*([C'_{d}]^\bullet)$, resp.
 $H^*([C_{d}]^\bullet)$. \nbox
\end{lemma}

\section{Freeness results}\label{s11}
The aim of this section is to establish some freeness proprieties for the homology bialgebras of $\TD$, $\TstD$ and of their duals.

\begin{theorem}\label{t111}
For any ring of coefficients divided powers operations are homology operations for complexes $\TD$, $\TstD$. \nbox  
\end{theorem}

\noindent {\bf Proof:} The main difficulty is to prove that divided powers of a border are also borders. Recall that $\TD$ and $\TstD$ are normalized parts of the cohomology
of cosimplicial spaces, see previous section. It means that these complexes are normalized parts of simplicial algebras. For any simplicial algebra one can define divided powers, 
and it is easy to see that our devided powers coincide with the standard divided powers defined for 
simplicial algebras (this is analogous to Lemma~\ref{l101}).
It is well-known, that divided powers are homotopy operations in this situation. In the case of the bar-construction this result is due to H.~Cartan~\cite{Cartan}. The proof is not more difficult
in the case of arbitrary simplicial algebra. Generalization for any type of simplicial algebras was obtained by B.~Fresse, \cf~\cite{Fresse1,Fresse2}, see Corollaries~2.2.11 and~2.2.12 
of~\cite{Fresse2}\footnote{Higher divided powers operations in the homotopy of simplicial algebras and relations between them are given in~\cite{Dwyer,GL} 
(over $\Z_{2}$) and in~\cite{Bousfield} (over any $\Z_{p}$).}.
\nbox

\begin{corollary}\label{c112}
The homology of $\TD$ and of $\TstD$ are Hopf algebras with divided powers. \nbox
\end{corollary}

\begin{theorem}\label{t113}
For any field of coefficients, the homology bialgebras of $\TD$, $\TstD$ are free Hopf algebras with divided powers. \nbox
\end{theorem}

\noindent {\bf Proof:} Recall that \lq\lq free Hopf algebra $\Gamma(V)$ with divided powers" generated by a graded vector space $V$ is the dual of
the polynomial bialgebra generated by the graded dual of $V$. It can be easily seen, that $\Gamma(V)$ can be endowed with divided powers, \cf~\cite{Andre}.
By Corollaries~\ref{c112} and~\ref{c54} the homology space of $\TD$ and that of $\TstD$ form graded bicommutative connected bialgebras with divided powers.
M.~Andr\'e proved that any connected Hopf algebra with divided powers\footnote{Existence of divided powers operations presumes that multiplicative structure is 
graded commutative.} is isomorphic to the universal enveloping coalgebra of a graded Lie coalgebra, \cf~\cite{Andre}\footnote{See also~\cite{Sj, Bl, Patras} for 
similar results.}. On the other hand, the cobracket of the corresponding
Lie coalgebras must be trivial, since our homology bialgebras are cocommutative, see Remark-Corollary~\ref{c54}. This implies the result. \nbox

\begin{corollary}\label{c114}
For any field of coefficients homology bialgebras of the duals to $\TstD$, $\TD$ are polynomial. \nbox
\end{corollary}

This corollary means that the Hochschild homology bialgebras of the operads of Poisson algebras, Gerstenhaber algebras and Batalin-Vylkovisky algebras are polynomial.

\begin{remark}\label{r114}
{\rm Note that in the case of characteristic 2, devided powers were defined for any elements of $\TstD$, $\TdstD$ both of even and of odd degrees. It means that in this case, 
the homology bialgebras of the duals to $\TstD$, $\TD$ are polynomial in the sense that odd degree generators are also polynomial (and not exterior). \nbox} 
\end{remark}

\begin{corollary}\label{c115}
For any field of coefficients bialgebras $\BB$ and $\BB_{0}$ of chord diagrams and their super-analogues $\tilde\BB$, $\tilde\BB_{0}$ are polynomial. \nbox
\end{corollary}

Recall that $\BB$, $\BB_{0}$, $\tilde\BB$, $\tilde \BB_{0}$ are lower-diagonal subbialgebras of the homology of the duals to $\TD^{odd}$, resp. $\TstD^{odd}$,
resp. $\TD^{even}$, resp. $\TstD^{even}$, see Section~\ref{s9}.

A similar result was obtained by S.~K.~Lando, \cf~\cite{Lando}, which asserts that for any ring of coefficents bialgebras $\BB$ and $\BB_{0}$ are generated
by their primitive elements.

\end{document}